\documentstyle[11pt, pb-diagram]{article}
\input amssym.def
\title{A geometric estimate for a periodic Schr\"odinger operator whose potential is the curvature of a spherical curve. \footnote {Supported by the SFB 288 of the DFG.}}
\author{THOMAS FRIEDRICH (Berlin)\\
}
\date{\today}

\begin{document}

\maketitle

\mbox{} \hrulefill \mbox{}\\

\newcommand{\vol}{\mbox{vol} \, }
\newcommand{\grad}{\mbox{grad} \, }

\begin{abstract}
We estimate from below by geometric data the eigenvalues of the periodic Sturm-Liouville operator $\displaystyle - 4 \frac{d^2}{ds^2} + \kappa^2 (s)$ with potential given by the curvature of a closed curve.
\end{abstract}

\vspace{0.5cm}

{\small
{\it Subj. Class.:} Differential geometry.\\
{\it 1991 MSC:} 58G25, 53A05.\\
{\it Keywords:} Dirac operator, spectrum, surfaces, Schr\"odinger operators, Fenchel inequality.} \\

\mbox{} \hrulefill \mbox{}\\
\section{Introduction}

\newfont{\graf}{eufm10}
\newcommand{\alth}{\mbox{\graf h}}

Let $X^3 (c)$ be a 3-dimensional space form of constant curvature $c=0$ or 1 and admitting a real Killing spinor with respect to some spin structure. Consider a compact, oriented and immersed surface $M^2 \subset X^3(c)$ with mean curvature $H$. The spin structure of $X^3 (c)$ induces a spin structure on $M^2$. Denote by $D$ the corresponding Dirac operator acting on spinor fields defined over the surface $M^2$. The first eigenvalue $\lambda_1^2 (D)$ of the operator $D^2$ and the first eigenvalue $\mu_1$ of the Schr\"odinger operator $\Delta + H^2 +c$ are related by the inequality

\[ \lambda_1^2 (D) \le \mu_1 (\Delta + H^2 + c) . \]

Equality holds if and only if the mean curvature $H$ is constant (see [1], [5]). Moreover, the Killing spinor defines a map $f \mapsto \Phi (f)$ of the space $L^2 (M^2)$ of functions into the space $L^2 (M^2; S)$ of spinors such that

\[ ||D(\Phi (f)) ||^2_{L^2} = \langle \Delta f + H^2f+ c f,f \rangle_{L^2} . \]

In particular, the mentioned inequality holds for all eigenvalues, i.e., 

\[ \lambda_k^2 (D) \le \mu_k (\Delta + H + c) . \]

This inequality was used in order to estimate the first eigenvalue of the Dirac operator defined on special surfaces of Euclidean space (see [1]). On the other hand, in case we know $\lambda_1^2 (D)$, the inequality yields  a lower bound for the spectrum of the Schr\"odinger operator $\Delta + H^2 + c$. For example, for any Riemannian metric $g$ on the 2-dimensional sphere $S^2$ the inequality

\[ \lambda_1^2 (D) \ge \frac{4 \pi}{\vol (S^2,g)} \]

holds (see [2], [6]). Consequently, we obtain

\[ \frac{4 \pi}{\vol (M^2,g)} \le \mu_1 (\Delta + H^2) \]

for any surface $M^2 \hookrightarrow {\Bbb R}^3$ of genus zero in Euclidean space ${\Bbb R}^3$. In this note we expose the described idea and, in particular, we estimate the spectrum of special periodic Schr\"odinger operators where the potential is given by the curvature $\kappa$ of a spherical curve. \\

\section{The one-dimensional case}

First of all, let us consider the 1-dimensional case, i.e., a curve $\gamma$ of length $L$ in a two-dimensional space form $X^2 (c)$. Let $\Phi$ be a Killing spinor of length one on $X^2 (c)$:

\[ \nabla_{T} {\Phi} = \frac{1}{2} c \cdot T \cdot \Phi . \]

The restriction $\varphi = \Phi_{|\gamma}$ defines a pair of spinors on $\gamma$ and the spinor field $\psi= f \cdot \varphi$ satisfies the equation:

\[ |D \psi |^2 = | \dot{f}|^2 +f^2 \left( \frac{c}{4} + \frac{1}{4} \kappa^2_g \right) , \]

where $\kappa_g$ is the curvature of the curve $\gamma$ in $X^2(c)$. Therefore, we obtain

\[ \lambda_k^2 (D) \le  \mu_k \left( - \frac{d^2}{ds^2} + \frac{c}{4}  + \frac{1}{4} \kappa^2_g \right) . \]

Suppose now that the spin structure on $\gamma$ induced by the spin structure of $X^2 (c)$ is non-trivial. Then we have $\displaystyle \lambda_{k+1}^2 (D)= \frac{4 \pi^2}{L^2} (k + 1/2)^2$ (see [4]) and, in particular,  we obtain

\[ \frac{4 \pi^2}{L^2} \left( k + \frac{1}{2} \right)^2 \le \mu_{k+1} \left( - \frac{d^2}{ds^2} + \frac{c}{4} + \frac{1}{4} \kappa^2_g  \right) . \]

{\bf Theorem 1:} {\it Let $\gamma \subset {\Bbb R}^3$ be a plane or spherical curve and denote by $\kappa^2 = c + \kappa_g^2$ the square of its curvature. Suppose that the induced spin structure on $\gamma$ is non-trivial, i.e., the tangent vector field has an odd rotation number. Then the inequality

\[ \frac{4 \pi^2}{L^2} \le \mu_1 \left( - 4 \frac{d^2}{ds^2} + \kappa^2 \right) \]

holds, where $\mu_1$ is the first eigenvalue of the periodic Sturm-Liouville operator on the interval $[0,L]$. Moreover, equality occurs if and only if the curvature is constant.}\\

{\bf Remark:} No geometric lower bound for the Sturm-Liouville operator $ - 4 \frac{d^2}{ds^2} + \kappa^2$ with potential defined by the square of the curvature $\kappa (s)$ of a closed curve $\gamma$ in Euclidean space seems to be known. We conjecture that the estimate given in \mbox{Theorem 1} holds for any closed curve in ${\Bbb R}^3$. Let us  compare this inequality with the well-known Fenchel-Milnor inequality

\[ 2 \pi \le \oint\limits_{\gamma} \kappa . \]

Thus, by the Cauchy-Schwarz inequality we obtain

\[ \frac{4 \pi^2}{L^2} \le  \frac{1}{L} \oint\limits_{\gamma} \kappa^2   . \]

Moreover, using the test function $f \equiv 1$, we have

\[ \mu_1 \left( - 4 \frac{d^2}{ds^2} + \kappa^2 \right) \le \frac{1}{L} \oint\limits_{\gamma} \kappa^2  . \]

Suppose that $\gamma$ is a simple curve in ${\Bbb R}^3$ and denote by $ \rho$ the minimal  number of generators of the fundamental group $\pi_1 ({\Bbb R}^3 \backslash \gamma)$. Then we have

\[ 2 \pi \rho \le \oint\limits_{\gamma} \kappa . \]

In the spirit of this remark one should be able to prove the stronger inequality 

\[ \frac{4 \pi^2}{L^2} \rho^2 \le \mu_1 \left( - 4 \frac{d^2}{ds^2} + \kappa^2 \right) \]

in case of a simple curve in ${\Bbb R}^3$. \\

{\bf Examples:} We calculated the eigenvalue $\mu_1$ for some classical curves in ${\Bbb R}^3$:

\begin{itemize}
\item[a.)] {\it The lemniscate:} $x= \sin (t) , \quad  y = \cos (t) \sin (t).$\\
\mbox{} \hspace{2.7cm} $4 \pi^2 /L^2 = 1.06193, \quad  \mu_1 = 3.7315, \quad \displaystyle \frac{1}{L} \oint\limits_{\gamma} \kappa^2  = 4.36004.$\\

\item[b.)] {\it The trefoil:} $x= \sin (3t)  \cos (t) , \quad  y = \sin (3t) \sin (t).$\\
\mbox{} \hspace{1.8cm} $4 \pi^2 /L^2 =0.221, \quad  \mu_1 = 5.21, \quad  \displaystyle \frac{1}{L} \oint\limits_{\gamma} \kappa^2 = 8.16.$\\

\item[c.)] {\it Viviani's curve:} $x= 1+ \cos (t) , \quad  y = \sin (2t) , \quad z= 2 \sin (t) .$\\
\mbox{} \hspace{2.6cm} $4 \pi^2 /L^2 =0.169071, \quad  \mu_1 = 0.5335, \quad  \displaystyle \frac{1}{L} \oint\limits_{\gamma} \kappa^2 = 0.567803.$\\

\item[d.)] {\it Torus knot:} $x= (8+3 \cos (5t)) \cos (2t) , \quad  y = (8+ 3 \cos (5t)) \sin (2t), \\ 
\mbox{} \hspace{2.3cm} z= 5 \sin (5t).$\\
\mbox{} \hspace{2.0cm} $4 \pi^2 /L^2 =0.00146034, \quad  \mu_1 = 0.03232, \quad  \displaystyle \frac{1}{L} \oint\limits_{\gamma} \kappa^2 = 0.0333803.$\\

\item[e.)] {\it The spherical spiral:} $x= \cos (t) \cos (4t) , \quad  y = \cos (t) \sin(4t), \quad z= \sin(t).$\\
\mbox{} \hspace{3.4cm} $4 \pi^2 /L^2 =0.127036, \quad  \mu_1 = 1.744, \quad  \displaystyle \frac{1}{L} \oint\limits_{\gamma} \kappa^2 = 4.93147.$\\

\end{itemize}

\section{The two-dimensional Schr\"odinger operator}

We generalize this inequality to the case of the two-dimensional periodic Schr\"odinger operator

\[ P_{A,L} = - \left( 1+ \frac{A^2}{L^2} \right) \frac{\partial}{\partial t^2} - 4 \frac{\partial^2}{\partial s^2} - \frac{4 A}{L} \frac{\partial}{\partial t} \frac{\partial}{\partial s}  + \kappa^2 (s) \]

defined on $[0,2\pi] \times [0,L]$:\\

{\bf Theorem 2:} {\it Let $\gamma \subset S^2 \subset {\Bbb R}^3$ be a closed, simple curve of length $L$ bounding a region of area $A$, and denote by $\kappa$ its curvature. Then the spectrum of the two-dimensional periodic Schr\"odinger operator $P_{A,L}$ is bounded by }

\[ \frac{4 \pi^2}{L^2} \le \mu_1 (P_{A,L}) . \]

{\it Equality holds if and only if the curvature of $\gamma$ is constant.}\\

In general, let us consider a Riemannian manifold $(Y^n,g)$ of dimension $n$ as well as an $S^1$-principal fibre bundle $\pi : P \to Y^n$ over $Y^n$. Denote by $\vec{V}$ the vertical vector field on $P$ induced by the action of the group $S^1$ on the total space $P$, i.e.,

\[ \vec{V} (p) = \frac{d}{dt} \left( p \cdot e^{it} \right)_{t=0} \quad , \quad p \in P . \]

A connection $Z$ in the bundle $P$ defines a decomposition of the tangent bundle $T(P)=T^v (P) \oplus T^h (P)$ into its vertical and horizontal subspace. We introduce a Riemannian metric $g^*$ on the total space $P$, requiring that 

\begin{itemize}
\item[a)] $g^* (\vec{V} , \vec{V})=1$,\\
\item[b)] $ g^* (T^v, T^h)=0$,\\
\item[c)] the differential $d \pi$ maps $T^h (P)$ isometrically onto $T(Y^n)$.
\end{itemize}

\bigskip

A closed curve $\gamma : [0,L] \to Y^n$ of length $L$ defines a torus $H(\gamma) := \pi^{-1} (\gamma) \subset P$ and we want to study the isometry class of this flat torus in $P$. Let $\alpha = e^{i \Theta} \in S^1 $ be the holonomy of the connection $Z$ along the closed curve $\gamma$. Consider a horizontal lift $\hat{\gamma} :[0,L] \to P$ of the curve $\gamma$. Then

\[ \hat{\gamma} (L) = \hat{\gamma} (0) e^{i \Theta} \]

holds. Consequently, the formula

\[ \Phi (t,s) = \hat{\gamma} (s) e^{-i \Theta s / L } e^{it} \]

defines a parametrization $\Phi : [0,2 \pi] \times [0,L] \to H(\gamma)$ of the torus $H(\gamma)$. Since

\[ \frac{\partial \Phi}{\partial t} = \vec{V} \quad , \quad \frac{\partial \Phi}{\partial s} = dR_{e^{it} e^{-i \Theta s/ L}} (\dot{\hat{\gamma}} (s)) - \frac{\Theta}{L} \vec{V} , \]

we obtain

\[ g^* \left( \frac{\partial \phi}{\partial t} , \frac{\partial \phi}{\partial t} \right) = 1  , \quad g^* \left( \frac{\partial \phi}{\partial t} , \frac{\partial \phi}{\partial s} \right) = - \frac{\Theta}{L} , \quad g^* \left( \frac{\partial \phi}{\partial s} , \frac{\partial \phi}{\partial s} \right) = 1 + \frac{\Theta^2}{L^2} , \]

i.e., the torus $H(\gamma)$ is isometric to the flat torus $({\Bbb R}^2 / \Gamma_o , g^*)$, where $\Gamma_o$ is the orthogonal lattice $\Gamma_o = 2 \pi \cdot {\Bbb Z} \oplus L \cdot {\Bbb Z}$ and the metric $g^*$ has the non-diagonal form

\[ g^* = \left( \begin{array}{cc} 1 & - \frac{\Theta}{L} \\  \\- \frac{\Theta}{L} & 1+ \frac{\Theta^2}{L^2} \end{array} \right) . \]

Using the transformation

\[ x = - \frac{\Theta}{L} s + t \quad , \quad y=s ,  \]

we see that $H(\gamma)$ is isometric to the flat torus $({\Bbb R}^2 / \Gamma , \, \, \, dx^2 + dy^2 )$,  where the lattice $\Gamma$ is generated by the two vectors

\[ v_1 = \left( \begin{array}{c} 2 \pi\\0 \end{array} \right) \quad , \quad v_2 = \left( \begin{array}{c} \Theta \\ L \end{array} \right) . \]

In case the closed curve $\gamma : [0,L] \to Y^n$ is the oriented boundary of an oriented compact surface $M^2 \subset Y^n$, we can calculate the holonomy $\alpha = e^{i \Theta}$ along the curve $\gamma$. Indeed, let $\Omega^Z$ be the curvature form of the connection $Z$. $\Omega^Z$ is a 2-form defined on the manifold $Y^n$ with values in the Lie algebra of the group $S^1$, i.e., with values in $i \cdot {\Bbb R}^1$. The parameter $\Theta$ is given by the integral 

\[ \Theta = i \int\limits_{M^2} \Omega^Z  . \]

Let us consider the Hopf fibration $\pi: S^3 \to S^2$,  where 

\[ S^3 = \{ (z_1, z_2) \in {\Bbb C}^2 : |z_1|^2 + |z_2|^2 =1 \} \]

is the 3-dimensional sphere of radius 1. The connection $Z$ is given by the formula

\[ Z= \frac{1}{2} \{ \bar{z}_1 dz_1  - z_1 d \bar{z}_1 + \bar{z}_2 dz_2 - z_2 d \bar{z}_2 \} \]

and its curvature form ($\omega = z_1 / z_2$)

\[ \Omega^Z = - \frac{d \omega \wedge d \bar{\omega}}{(1+|\omega|^2)^2} = - \frac{i}{2} dS^2 \]

essentially coincides with one half of the volume form of the unit sphere $S^2$ of radius 1. However, the differential $d \pi : T^h (S^3) \to T(S^2)$ multiplies the length of a vector by two, i.e., the Hopf fibration is a Riemannian submersion in the sense described before if we fix  the metric of the sphere $S^2 (\frac{1}{2})= \{ x \in {\Bbb R}^3 : |x| = \frac{1}{2} \}$ on $S^2$. Consequently, in case of a closed simple curve $\gamma \subset S^2$  bounding a region of area $A$,  the Hopf torus $H (\gamma) \subset S^3$ is isometric to the flat torus ${\Bbb R}^2 / \Gamma$ and the lattice $\Gamma$ is generated by the two vectors

\[ v_1 = \left( \begin{array}{c} 2 \pi\\0 \end{array} \right) \quad , \quad v_2 = \left( \begin{array}{c}  {A}/{2} \\  {L}/{2} \end{array} \right) . \]

 The mean curvature $H$ of the torus $H(\gamma) \subset S^3$ coincides with the geodesic curvature $\kappa_g$ of the curve $\gamma \subset S^2 \subset {\Bbb R}^3$ (see [7], [8]). We apply now the inequality

\[ \lambda_1^2 (D) \le \mu_1 ( \Delta + H^2 + 1) \]

to the Hopf torus $H (\gamma) \subset S^3$. Then we obtain the estimate

\[ \lambda_1^2 (D) \le \mu_1 \left( P_{A,L}  \right) , \]

where $D$ is the Dirac operator on the flat torus ${\Bbb R}^2 / \Gamma$ with respect to the induced spin structure. All spin structures of a 2-dimensional torus are classified by pairs $(\varepsilon_1, \varepsilon_2)$ of numbers $\varepsilon_i =0,1$. If $\gamma$ is a simple curve in $S^2$,  the induced spin structure on the Hopf torus $H(\gamma)$ is non-trivial and given by the pairs $(\varepsilon_1, \varepsilon_2) =(0,1)$. The spectrum of the Dirac operator for all flat tori is well-known (see [4]): The dual lattice $\Gamma^*$ is generated by

\[ v_1^* = \left( \begin{array}{c} \displaystyle \frac{1}{2 \pi} \\ \\ \displaystyle  - \frac{A}{2 \pi L}  \end{array} \right) \quad , \quad v_2^* = \left( \begin{array}{c} 0 \\ \\ \displaystyle \frac{2}{L} \end{array} \right)  \]

and the eigenvalues of $D^2$ are given by

\begin{eqnarray*}
\lambda^2 (k,l) &=& 4 \pi^2 \Big| \Big| \, \, k  v_1^* + \left(l + \frac{1}{2}\right) v_2^* \Big| \Big|^2 =\\
&=& k^2 + \frac{4 \pi^2}{L^2} \left((2l+1) - k \frac{A}{2 \pi} \right)^2 .
\end{eqnarray*}

We minimize $\lambda^2 (k,l)$ on the integral lattice ${\Bbb Z}^2$.  The isoperimetric inequality $4 \pi A - A^2 \le L^2$ and $A \le \vol (S^2) = 4 \pi$ yield the result that $\lambda^2 (k,l)$ attends its minimum at $(k,l) = (0,1)$, i.e., 

\[  \frac{4 \pi^2}{L^2} \le \lambda^2 (k,l)  . \]

{\bf Remark 1:} We replace the Hopf fibration by the $S^1$-principal fibre bundle of Chern class $m \ge 0$. The corresponding total space is the Lens space $L(m,1)$ and we have the commutative diagram

\[
\begin{diagram}
\node{S^3} \arrow{se,b}{\pi}  \arrow[2]{e} \node{} \node{L(m,1)} \arrow{sw,b}{\pi_m}\\
\node{} \node{S^2} \node{}
\end{diagram}
\]

Let $H_m (\gamma) \subset L(m,1)$ be the Hopf torus. $H_m (\gamma)$ is isometric to ${\Bbb R}^2 / \Gamma_m$, where the lattice $\Gamma_m$ is generated by the vectors

\[ v_1 = \left( \begin{array}{c} 2 \pi / m \\ 0 \end{array} \right) \quad , \quad v_2 = \left( \begin{array}{c} A/2 \\ L/2 \end{array} \right) . \]

Moreover, the Lens space $L(m,1)$ admits a unique spin structure with a Killing spinor (see [3]). Even in case of $m \not= 1$,  the induced spin structure on $H_m (\gamma)$ is described by the parameters $(\varepsilon_1 , \varepsilon_2)=(0,1)$. Since the local geometry of $H_m (\gamma)$ in $L(m,1)$ essentially coincides with the geometry of $H(\gamma)$ in $S^3$, we obtain the inequality

\[ \frac{4 \pi^2}{L^2} = \min\limits_{(k,l)} \left\{ k^2m^2 + \frac{4 \pi^2}{L^2} \left( (2l+1) - k \frac{mA}{2 \pi} \right)^2 \right\} \le \mu_1 \left( - 4 \frac{d^2}{ds^2} + \kappa^2 \right) . \]

Consequently, the investigation of the two-dimensional Schr\"odinger operator in case of $m \not= 1$ yields the same result for the Sturm-Liouville operator as above. \\

{\bf Remark 2:} Suppose now that  equality holds for some curve $\gamma \subset S^2$. We consider the corres\-ponding  Hopf torus $H(\gamma) \subset S^3$ and then we obtain

\[ \lambda_1^2 (D) = \mu_1 ( \Delta + H^2 +1) . \]

Therefore, the mean curvature $H= \kappa$ is constant, i.e., $\gamma$ is a curve on $S^2$ of constant curvature $\kappa$. Consequently, $\gamma$ is a circle in a 2-dimensional plane. Denote by $r$ its radius. Then

\[ \kappa^2 = \frac{1}{r^2} \quad , \quad L= 2 \pi r \quad , \quad A=2 \pi (1 - \sqrt{1-r^2})  , \]

and the inequality

\[ \frac{4 \pi^2}{L^2} \le \kappa^2   \]

is an  equality for all $r \not= 0$. \\


\bigskip


\small
THOMAS FRIEDRICH\\
Humboldt-Universit\"at zu Berlin, Institut f\"ur Mathematik, Sitz: Ziegelstra\ss e 13a, \\
Unter den Linden 6, D-10099 Berlin\\
{\tt e-mail: friedric@mathematik.hu-berlin.de}\\

\end{document}